\documentclass[12pt, reqno]{amsart}
\usepackage{amsmath, amsthm, amscd, amsfonts, amssymb, graphicx, color}
\usepackage[bookmarksnumbered, colorlinks, plainpages]{hyperref}

\newcommand{\R}{\mbox{$\mathbb{R}$}}
\newcommand{\N}{\mbox{$\mathbb{N}$}}

\begin{document}
\thispagestyle{empty}
 \setcounter{page}{1}

\begin{center}
{\large\bf An extension  of  Mizoguchi--Takahaashi's fixed point
theorem

\vskip.15in

{\bf  M. Eshaghi Gordji, H. Baghani, M. Ramezani and H. Khodaei }
 \\[2mm]

{\footnotesize Department of Mathematics, Semnan University,\\ P. O.
Box 35195-363, Semnan, Iran}}
\end{center}
\vskip 5mm

 \noindent{\footnotesize{\bf Abstract.}
 Our main theorem
is an extension of the well--known Mizoguchi--Takahaashi's fixed
point theorem [N. Mizogochi and W. Takahashi, Fixed point theorems
for multi--valued mappings on complete metric space,
 {\it J. Math. Anal. Appl.} 141 (1989) 177--188].

 \vskip.15in

\footnotetext{\textbf{E-mail}{\tt :  madjid.eshaghi@gmail.com,
h.baghani@gmail.com, ramezanimaryam873@gmail.com,
khodaei.hamid.math@gmail.com}}

}

  \newtheorem{df}{Definition}[section]
  \newtheorem{rk}[df]{Remark}
   \newtheorem{lem}[df]{Lemma}
   \newtheorem{thm}[df]{Theorem}
   \newtheorem{pro}[df]{Proposition}
   \newtheorem{cor}[df]{Corollary}
   \newtheorem{ex}[df]{Example}

 \setcounter{section}{0}
 \numberwithin{equation}{section}

\vskip .2in
\begin{center}
\section{INTRODUCTION AND STATEMENT OF RESULTS}
\end{center}

Let $(X,d)$ be a  metric space. $CB(X)$ denotes the collection of
all nonempty closed bounded subsets of $X$. For $A,B \in CB(X)$,
and $x\in X$, define $D(x,A):=\inf\{d(x,a);a\in A\}$, and
$$H(A,B):=\max\{\sup_{a\in A}D(a,B),\sup_{b\in B}D(b,A).$$ It is easy to see that
$H$ is a metric on $CB(X)$. $H$ is called the Hausdorff metric
induced by $d$.

\begin{df} An element $x\in X$ is said to be a fixed point of a multi--valued mapping
$T:X\rightarrow CB(X)$, if such that $x\in T(x)$.
\end{df}

One can show that $(CB(X),H)$ is a complete metric space,
whenever $(X,d)$ is a complete metric space (see for example
Lemma $8.1.4$, of \cite{Rus}).

In 1969, Nadler  \cite {N} extended the Banach contraction principle
 \cite {Ba} to set--valued mappings as follows.

\begin{thm}   Let $(X, d)$
be a complete metric space and let $T$ be a mapping from $X$ into
$CB(X)$. Assume that there exists $r\in [0, 1)$ such that $\mathcal
H_d(T x,Ty)\leq rd(x,y)$ for all $x,y\in X$. Then there exists $z\in
X$ such that $z\in T(z).$
\end{thm}

Nadler's theorem  was generalized by Mizoguchi and Takahaashi
\cite{A} in the following way.

\begin{thm}
Let $(X,d)$ be a complete metric space and  let $T$ be a mapping
from $(X,d)$ into $(CB(X), H)$ satisfies
$$H(Tx,Ty)\leq\alpha(d(x,y))d(x,y)$$ for all $x,y \in X$, where
$\alpha$ be a function from $[0,\infty)$ into $[0,1)$ such that
$\limsup_{s\to t^+}\alpha(s)<1$ for all $t\in [0,\infty)$. Then $T$
has a fixed point.
\end{thm}
Recently Suzuki \cite{Su} proved the Mizoguchi--Takahashi's fixed
point theorem by an interesting and short proof.\\

On the other hand, Banach contraction principle was generalized by
Reich  \cite{R1,R2} as follows.

\begin{thm}
Let $(X,d)$ be a complete metric space and  let $T$ be a mapping
from $(X,d)$ into $(CB(X), H)$ satisfies $$H(Tx,Ty)\leq
\beta[D(x,Tx)+D(y,Ty)]$$ for all $x,y \in X$, where
$\beta\in[0,\frac{1}{2})$. Then $T$ has a fixed point.
\end{thm}

In 1973, Hardy and Rogers  \cite {HR} extended the Reich's theorem
by the following way.

\begin{thm}
Let $(X,d)$ be a complete metric space and let $T$ be a mapping from
$X$ into $X$ such that $$d(Tx,Ty)\leq \alpha d(x,y)+\beta
[d(x,Tx)+d(y,Ty)]+\gamma [d(x,Ty)+d(y,Tx)]$$ for all $x,y \in X$,
where $\alpha,\beta,\gamma \geq 0$ and $\alpha+2\beta+2\gamma <1$.
Then $T$ has a fixed point.
\end{thm}
Recently, the authors of the present paper \cite {EBKR} extended the
theorems 1.5 and 1.2 as follows.
\begin{thm}
Let $(X,d)$ be a complete metric space and  let $T$ be a mapping
from $X$ into $CB(X)$ such that $$H(Tx,Ty)\leq \alpha d(x,y)+\beta
[D(x,Tx)+D(y,Ty)]+\gamma [D(x,Ty)+D(y,Tx)]$$ for all $x,y \in X$,
where $\alpha,\beta,\gamma \geq 0$ and $\alpha+2\beta+2\gamma <1$.
Then $T$ has a fixed point.
\end{thm}

In this paper, we shall  generalize above results.  More precisely,
we prove the following theorem, which can be regarded as an
extension of all theorems 1.2,1.3,1.4,1.5 and 1.6.

\begin{thm}
Let $(X,d)$ be a complete metric space and  let $T$ be mapping
from $X$ into $CB(X)$ such that
\begin{eqnarray*}
H(Tx,Ty)&\leq& \alpha (d(x,y))d(x,y)+\beta(d(x,y))
[D(x,Tx)+D(y,Ty)]\\&&+\gamma(d(x,y)) [D(x,Ty)+D(y,Tx)]
\end{eqnarray*}
 for all $x,y \in X$, where $\alpha,\beta,\gamma$ are mappings from
$[0,\infty)$ into $[0,1)$ such that
$\alpha(t)+2\beta(t)+2\gamma(t)<1$ and $\limsup_{s\to
t^+}\frac{\alpha(t)+\beta(t)+\gamma(t)}{1-(\beta(t)+\gamma(t))}<1$
for all $t\in[0,\infty)$. Then $T$ has a  fixed point.
\end{thm}
Moreover, we conclude the following results by using  theorem 1.7.
\begin{cor}
Let $(X,d)$ be a complete metric space and  let $T$ be a mapping
from $(X,d)$ into $(CB(X), H)$ satisfies $$H(Tx,Ty)\leq
\beta(d(x,y))[D(x,Tx)+D(y,Ty)]$$ for all $x,y \in X$, where $\beta$
be a function from $[0,\infty)$ into $[0,\frac{1}{2})$ and
$\limsup_{s\to t}\beta(s)<\frac{1}{2}$ for all $t\in [0,\infty)$.
Then $T$ has a fixed point.
\end{cor}

\begin{cor}
Let $(X,d)$ be a complete metric space and  let $T$ be a mapping
from $(X,d)$ into $(CB(X), H)$ satisfies
$$H(Tx,Ty)\leq\alpha(d(x,y))d(x,y)+ \beta(d(x,y))[D(x,Tx)+D(y,Ty)]$$ for
all $x,y \in X$, where $\alpha,\beta$ are function from $[0,\infty)$
into $[0,1)$ such that $\alpha(t)+2\beta(t)<1$ and $\limsup_{s\to
t^+}(\frac{\alpha(t)+\beta(t)}{1-\beta(t)})<1$ for all $t\in
[0,\infty)$. Then $T$ has a fixed point.
\end{cor}

\vskip .2in
\begin{center}
\section{Proof of the main theorem}
\end{center}

\begin{proof}
Define function $\alpha^{'}$ from $[0,\infty)$ into $[0,1)$ by
$\alpha^{'}(t)=\frac{\alpha(t)+1-2\beta(t)-2\gamma(t)}{2}$ for
$t\in [0,\infty)$. Then we have the following assertions:\\\\
$1)$ $\alpha(t)<\alpha^{'}(t)$ for all $t\in[0,\infty)$.\\\\
$2)$ $\limsup_{s\to
t^+}\frac{\alpha^{'}(t)+\beta(t)+\gamma(t)}{1-(\beta(t)+\gamma(t))}<1$ for all $t\in[0,\infty)$.\\\\
$3)$ For $x,y\in X$ and $u\in Tx$, there exists $\nu\in Ty $ such
that
\begin{eqnarray*}
d(\nu,u)&\leq&\alpha^{'} (d(x,y))d(x,y)+\beta(d(x,y))
[D(x,Tx)+D(y,Ty)]\\&&+\gamma(d(x,y)) [D(x,Ty)+D(y,Tx)].
\end{eqnarray*}
Putting $u=y$ in 3), we obtain that:\\
$4)$ For $x\in X $ and $y\in Tx$ there exists $\nu\in Ty$ such that
\begin{eqnarray*}
d(\nu,y)&\leq& \alpha^{'} (d(x,y))d(x,y)+\beta(d(x,y))
[D(x,Tx)+D(y,Ty)]\\&&+\gamma(d(x,y)) [D(x,Ty)+D(y,Tx)].
\end{eqnarray*}
Hence,  we can define sequence $\{x_n\}_{n\in\N}$ such that
$x_{n+1}\in Tx_n, x_{n+1}\neq x_n$ and
\begin{eqnarray*}
d(x_{n+2},x_{n+1})&\leq& \alpha^{'}
(d(x_{n+1},x_{n}))d(x_{n+1},x_{n})+\beta(d(x_{n+1},x_{n}))
[D(x_n,Tx_n)\\&&+D(x_{n+1},Tx_{n+1})]+\gamma(d(x_{n+1},x_{n})
[D(x_n,Tx_{n+1})\\&&+D(x_{n+1},Tx_n)]
\end{eqnarray*}for all $n\in\N$. It follows that

$$
d(x_{n+2},x_{n+1})\leq\frac{\alpha^{'}(d(x_{n+1},x_{n}))+\beta(d(x_{n+1},x_{n}))+\gamma(d(x_{n+1},x_{n}))}
{1-(\beta(d(x_{n+1},x_{n}))+\gamma(d(x_{n+1},x_{n})))}d(x_{n+1},x_{n})$$for
all $n\in\N$. On the other hand, we have
$$\frac{\alpha^{'}(t)+\beta(t)+\gamma(t)}{1-(\beta(t)+\gamma(t))}<1$$
for all $t\in\ [0,\infty)$, then $\{d(x_{n+1},x_n)\}$ is a
non-increasing sequence in $\R$. Hence, $\{d(x_{n+1},x_n)\}$ is a
converges to some nonnegative integer $\tau$. By assumption,
$$\limsup_{s\to \tau^+
}\frac{\alpha^{'}(s)+\beta(s)+\gamma(s)}{1-(\beta(s)+\gamma(s))}<1$$
so, we have
$$\frac{\alpha^{'}(\tau)+\beta(\tau)+\gamma(\tau)}{1-(\beta(\tau)+\gamma(\tau))}<1$$
then, there exist $r\in [0,1)$ and $\epsilon>0$ such that
$$\frac{\alpha^{'}(s)+\beta(s)+\gamma(s)}{1-\beta(s)+\gamma(s)}<r$$
for all $s\in[\tau,\tau+\epsilon]$. We can take $\nu\in\N$ such that
$$\tau\leq d(x_{n+1},x_n)\leq\tau+\epsilon$$ for all $n\in\N$ with
$n\geq\nu$. It follows that
\begin{eqnarray*}
d(x_{n+2},x_{n+1})&\leq&\frac{\alpha^{'}(d(x_{n+1},x_{n}))+\beta(d(x_{n+1},x_{n}))+\gamma(d(x_{n+1},x_{n}))}{1-(\beta(d(x_{n+1},x_{n}))+\gamma(d(x_{n+1},x_{n})))}d(x_{n+1},x_{n})\\&\leq&
r d(x_{n+1},x_{n})
\end{eqnarray*}
 for all $n\in\N$ with $n\geq\nu$. This implies  that \\
\begin{eqnarray*}
\sum_{n=1}^\infty d(x_{n+2},x_{n+1})\leq\sum_{n=1}^\nu
d(x_{n+1},x_n)+\sum_{n=1}^\infty r^n d(x_{\nu+1},x_{\nu})<\infty.
\end{eqnarray*}
Hence, $\{x_n\}$ is a Cauchy sequence. Since $X$ is a complete
metric space, then $\{x_n\}$  converges to some point $x^*\in X$.
Now, we have
\begin{eqnarray*}
D(x^*,Tx^*)&\leq& d(x^*,x_{n+1})+D(x_{n+1},Tx^*)\\&\leq&
d(x^*,x_{n+1})+H(Tx_{n},Tx^*)\\&\leq&
d(x^*,x_{n+1})+\alpha(d(x_{n},x^*))
d(x_{n},x^*)\\&+&\beta(d(x_{n},x^*))
[D(x_{n},Tx_{n})+D(x^*,Tx^*)]\\&+&\gamma(d(x_{n},x^*))
[D(x_{n},Tx^*)+D(x^*,Tx_{n})]
\end{eqnarray*}
for all $n\in\N$. Therefore,
\begin{eqnarray*}
D(x^*,Tx^*)&\leq& d(x^*,x_{n+1})+\alpha(d(x_{n},x^*))
d(x_{n},x^*)\\&+&\beta(d(x_{n},x^*))
[d(x_{n+1},x_{n})+D(x^*,Tx^*)]\\&+&\gamma(d(x_{n},x^*))
[D(x_{n},Tx^*)+d(x_{n},x^*)]
\end{eqnarray*}for all $n\in\N$.
It follows that
\begin{eqnarray*}
D(x^*,Tx^*)&\leq&\liminf_{n\to\infty}(\beta(d(x_n,x^*))+\gamma(d(x_n,x^*)))
D(x^*,Tx^*)\\&=&\liminf_{s\to 0^+}(\beta(s)+\gamma(s))
D(x^*,Tx^*)\\&\leq&\limsup_{s\to 0^+
}(\frac{\alpha(s)+\beta(s)+\gamma(s)}{1-(\beta(s)+\gamma(s))})
D(x^*,Tx^*).
\end{eqnarray*}
On the other hand, we have  $$\limsup_{s\to
0^+}(\frac{\alpha(s)+\beta(s)+\gamma(s)}{1-(\beta(s)+\gamma(s))})<1$$
then $D(x^*,Tx^*)=0$. Since $Tx^*$ is closed, then $x^*\in Tx^*$.
\end{proof}

\end{document}